\documentclass[12pt]{article}
\usepackage{amsthm,amsmath,amssymb}
\usepackage{fancyvrb}
\theoremstyle{plain}
\newtheorem{theorem}{Theorem}

\theoremstyle{definition}

\newcommand{\abs}[1]{\lvert#1\rvert}
\newcommand{\NN}{\mathbb N}
\DeclareMathOperator{\PGL}{PGL}
\DeclareMathOperator{\PSL}{PSL}
\usepackage{xcolor} 
\usepackage{url}
\usepackage[colorlinks,linkcolor=blue,citecolor=purple]{hyperref}
\begin{document}
\title{A note about solvable and non-solvable finite groups of the
  same order type}
\author{Peter M\"uller\\[1mm]
  Institute of Mathematics, University of W\"urzburg\\
  \texttt{peter.mueller@uni-wuerzburg.de}}
\maketitle
\begin{abstract}
  Two finite groups are said to have the same order type if for each
  positive integer $n$ both groups have the same number of elements of
  order $n$. In 1987 John G.~Thompson (see \cite[Problem
  12.37]{kourovka_notebook}) asked if in this case the solvability of
  one group implies the solvability of the other group.

  In 2024 Pawe\l{} Piwek gave a negative example in \cite{piwek2024solvablenonsolvablefinitegroups}. He
  constructed two groups of order
  $2^{365}\cdot3^{105}\cdot7^{104}\approx7.3\cdot10^{247}$ of the same
  order type, where only one is solvable.

  In this note we produce a much smaller example of order
  $2^{13}\cdot3^4\cdot7^3=227598336$.
\end{abstract}
\section{Introduction}
The order type of a finite group $G$ is the function
$o_G:\NN_{>0}\to\NN$ such that $o_G(n)$ is the number of elements of
order $n$ in $G$.

The following was an open problem since the late 80's, originally
raised by John G.~Thompson (see e.g.~\cite[Problem
12.37]{kourovka_notebook}): Let $G$ and $H$ be finite groups. If
$o_G=o_H$ and $G$ is solvable, is then $H$ necessarily solvable too?

In 2024 Pawe\l{} Piwek gave a negative example in \cite{piwek2024solvablenonsolvablefinitegroups}, where
$\abs{G}=\abs{H}=2^{365}\cdot3^{105}\cdot7^{104}\approx7.3\cdot10^{247}$.

We will give a much smaller counterexample. As to the notation: $C_i$
and $D_i$ denote the cyclic or dihedral group of order $i$,
respectively. $Q_i$ is the (generalized) quaternion group of order
$i$, and $A_4$ is the alternating group of degree $4$. The Id of a
group $S$ is the pair $(n,i)$ such that $S$ is isomorphic to the
$i$-th group of order $n$ in the list of small groups in
\cite{besche:small_groups}. In the computer algebra systems Gap
\cite{GAP4} and Magma \cite{magma}, this group is available as
\texttt{SmallGroup(n, i)}, or as \texttt{libgap.SmallGroup(n, i)} in
SageMath \cite{sagemath}. The following table lists seven groups
$S_1,\dots,S_7$ which we use to build the counterexample. Note that
the semidirect products in the structural description do not uniquely
determine the group up to isomorphism.
\begin{table}[h]
  \centering\[
  \begin{array}{ccc}
    &\text{Id}&\text{description}\\
    \hline
    S_1 & (168, 43) & C_2^3\rtimes(C_7\rtimes C_3)\\
    S_2 & (1008, 289) & C_7 \rtimes(C_3\times(C_3\rtimes Q_{16}))\\
    S_3 & (1344, 6967)& C_7\rtimes(((C_4\times D_8)\rtimes C_2)\rtimes C_3)\\
    S_4 & (21, 1) & C_7\rtimes C_3\\
    S_5 & (96, 166) & C_{12}\times Q_8\\
    S_6 & (336, 136) & C_7\rtimes(C_4\times A_4)\\
    S_7 & (336, 208) & \PGL_2(7)
  \end{array}\]
  \caption{Some small groups}
\end{table}

\begin{theorem}
  With $S_i$ as in the table, set $G=S_1\times S_2\times S_3$ and
  $H=S_4\times S_5\times S_6\times S_7$. Then $o_G=o_H$ and $G$ is
  solvable, but $H$ is not solvable.
\end{theorem}
\section{Proof of the theorem}
The proof follows the lines as in \cite{piwek2024solvablenonsolvablefinitegroups}. For a finite group $X$ and
$n\ge1$ let $e_X(n)$ be the number of $x\in X$ such with $x^n=1$. The
function $e_X$ is called the exponent type of $X$. Note that if $X$
and $Y$ are finite groups, then $o_X=o_Y$ if and only if $e_X=e_Y$
(see e.g.~\cite[Lemma 1]{piwek2024solvablenonsolvablefinitegroups}) and $e_{X\times Y}=e_X\cdot e_Y$ (see
e.g.~\cite[Lemma 2]{piwek2024solvablenonsolvablefinitegroups}).

Furthermore, $e_X(n)$ depends only on the residue of $n$ modulo the
exponent $m$ of $X$. In addition, $e_X(n)=e_X(\gcd(n,m))$. Thus $e_X$
is determined by $e_X(n)$ for the divisors $n$ of the exponent of $X$.

Clearly, $G$ as a direct product of solvable groups is
solvable. However, $H$ is not solvable because the derived subgroup
$S_7'=\PSL_2(7)$ is perfect.

Using a computer algebra system one verifies that $G$ and $H$ have
exponent $168$, and that
\[
  e_{S_1}(n)\cdot e_{S_2}(n)\cdot e_{S_3}(n)%
  =e_{S_4}(n)\cdot e_{S_5}(n)\cdot e_{S_6}(n)\cdot e_{S_7}(n)
\]
for each divisor $n$ of $168$, hence $e_G(n)=e_H(n)$ for these
$n$. The ancillary files \texttt{verification.sage},
\texttt{verification.gap} and \texttt{verification.mag} contain a
verification of this computation using SageMath, Gap, and Magma,
respectively. Note that for SageMath and Magma the online calculators
at \url{https://sagecell.sagemath.org/} and \url{http://magma.maths.usyd.edu.au/calc/} can be used.

As an example, here is the SageMath code for
\texttt{verification.sage}:
{\small
\begin{Verbatim}[samepage=true, frame=single]
def e(X, n):
    return [x^n == 1 for x in X.List()].count(True)

ids = [(168, 43), (1008, 289), (1344, 6967), (21, 1), (96, 166),
       (336, 136), (336, 208)]
    
l = [libgap.SmallGroup(n, i) for n, i in ids]
lg = l[:3] # The factors of G
lh = l[3:] # The factors of H

print(168 == lcm(int(x.Exponent()) for x in lg))
print(168 == lcm(int(x.Exponent()) for x in lh))

for d in divisors(168):
    print(prod(e(s, d) for s in lg) == prod(e(s, d) for s in lh))
  \end{Verbatim}
  }

\end{document}